\documentclass[10pt]{amsart}
\usepackage{mathptmx}
\usepackage{amsmath}
\usepackage{amssymb}
\usepackage{array}
\usepackage{geometry}
\usepackage[bookmarks=true,colorlinks=true, pdfstartview=FitV, linkcolor=black, citecolor=blue, urlcolor=black]{hyperref}

\usepackage{color}
\definecolor{DarkRed}{rgb}{0.55,.00,0.2}
\definecolor{DarkGrey}{rgb}{0.35,.35,0.35}

\theoremstyle{definition}

\theoremstyle{remark}

\numberwithin{equation}{section}

\begin{document}

\title{ On the generalized Lebedev  index transform}

\author{S. Yakubovich}
\address{Department of Mathematics, Faculty of Sciences,  University of Porto,  Campo Alegre str.,  687; 4169-007 Porto,  Portugal}
\email{ syakubov@fc.up.pt}

\keywords{Index Transform, Lebedev transform, Kontorovich-Lebedev transform, Macdonald function,  Meijer transform, Laplace transform, Fourier transform, Mellin transform, Initial value problem}
\subjclass[2000]{  44A15, 33C10, 44A05
}

\date{\today}
\maketitle

\markboth{\rm \centerline{ S.  Yakubovich}}{}
\markright{\rm \centerline{Generalized Lebedev Index Transform }}

\begin{abstract}  An essential generalization of the Lebedev index transform with the square of the Macdonald function is investigated. Namely, we consider a  family of integral operators with the positive kernel $\left|K_{(i\tau+\alpha)/2}(x)\right|^2,\\   \alpha \ge 0,\  x >0,   \  \tau   \in \mathbb{R},$ where $K_\mu(z)$ is the Macdonald function and $ i$ is the imaginary unit.  Mapping properties such as the boundedness, compactness, invertibility are investigated for these operators  and their adjoints  in the Lebesgue weighted spaces.  Inversion theorems are proved.   Important  particular cases are exhibited.   As an interesting application, a solution of the initial   value problem for the second order differential difference equation, involving the Laplacian,  is obtained. 
\end{abstract}

\section{Introduction and preliminary results}

Let $\alpha \ge 0$.  The main goal of this paper is to investigate mapping properties of a family of index transforms \cite{yak} and their adjoints, involving  the Macdonald function in  the kernel,  namely,
$$F_\alpha(\tau) = \int_0^\infty \left| K_{(i\tau+ \alpha)/2}(x)\right|^2  f(x)dx,   \quad    \tau \in \mathbb{R}, \eqno(1.1)$$
$$G_\alpha(x) = \int_{-\infty}^\infty \left| K_{(i\tau+ \alpha)/2}(x)\right|^2  g(\tau)d\tau,   \quad   x  \in \mathbb{R}_+, \eqno(1.2)$$
where $i$ is the imaginary unit.  The Macdonald function $K_\mu(z)$ \cite{erd}, Vol. II is the modified Bessel function of the second kind,   which satisfies the differential equation
$$  z^2{d^2u\over dz^2}  + z{du\over dz} - (z^2+\mu^2)u = 0.\eqno(1.3)$$
It has the asymptotic behaviour 
$$ K_\mu(z) = \left( \frac{\pi}{2z} \right)^{1/2} e^{-z} [1+
O(1/z)], \qquad z \to \infty,\eqno(1.4)$$
and near the origin
$$z^\mu K_\mu(z) = 2^{\mu-1} \Gamma(\mu) + o(1), \ z \to 0,\eqno(1.5)$$
$$K_0(z) = -\log z + O(1), \ z \to 0. \eqno(1.6)$$
The Macdonald function can be represented by the integral 
$$K_{\mu}(z)= \int_0^\infty e^{- z\cosh u}\cosh (\mu u)du,\   {\rm Re}\  z >0,  \   \mu \in \mathbb{C}.\eqno(1.7)$$
Concerning the product of the Macdonald  functions  $K_{(\mu+\alpha)/2}(z)K_{(\mu-\alpha)/2}(z)$ the key formula, which will be used in the sequel is relation (2.16.5.4) in \cite{prud},  Vol. II
$$K_{(\mu+\alpha)/2}(z)K_{(\mu-\alpha)/2}(z) = \int_0^\infty K_\mu\left(z\left( x + {1\over x}\right)\right) x^{\alpha-1} dx,
\quad   {\rm Re}\  z >0.\eqno(1.8)$$
Letting in (1.1), (1.2) $\alpha=0$, we come up with the operator of the Lebedev index transform and its adjoint, which is associated with the square of the Macdonald function \cite{square}, \cite{leb}.   We note that, indeed, an essential generalization of the Lebedev transform will be investigated, since it is impossible to reduce (1.1), (1.2) to the Lebedev operator  via any substitution of parameters or functions.   Other index transforms related to the product of Macdonald's functions of different arguments considered by the author in \cite{yaprod}, \cite{yaprodrev}.  Our method of investigation of the operators (1.1), (1.2) will involve a similar technique, which was employed to study the boundedness and   invertibility properties of the Kontorovich-Lebedev transform \cite{yakl} and  general index transforms \cite{gen}.  

\section {Boundedness and compactness in  Lebesgue's spaces }

Let us introduce the following Lebesgue functional spaces 
$$L^\alpha\equiv   L_1\left(\mathbb{R}_+; K^2_{\alpha/2}(x)dx\right) := 
\left\{f:  \int_0^\infty K^2_{\alpha/2}(x)|f(x)|dx < \infty \right\}.\eqno(2.1)$$
In particular, as we will show below,  it contains  spaces  $L_{\nu, p}(\mathbb{R}_+)$ for some $\nu \in \mathbb{R},\   1\le  p \le \infty$ with the norms
$$||f||_{\nu, p}=\left(\int_0^\infty x^{\nu  p-1}|f(x)|^pdx\right)^{1/p}< \infty,\eqno(2.2)$$
$$||f||_{\nu, \infty}=\hbox{ess sup}_{x\ge 0}   |x^\nu f(x)|< \infty.$$
When $\nu={1\over p}$ we obtain the usual norm in $L_p$ denoted by   $|| \ ||_p$.  

{\bf Lemma 1.} {\it Let  $\alpha \ge 0,\  \nu+\alpha <1, 1\le p \le \infty, \ q={p\over  p-1}$.   Then the embedding holds
$$L_{\nu, p}(\mathbb{R}_+)\subseteq L^\alpha \eqno(2.3)$$
and
$$ ||f||_{L^\alpha} \le \left[ {\Gamma^{1/q}\left(q(1-\nu)\right)\over
4 q^{1-\nu} } \  B\left({1-\nu\over 2}+ {\alpha\over 4},\  {1-\nu\over 2}- {\alpha\over 4}\right) \right]^2   ||f||_{\nu, p}, \ 1<p\le  \infty, \eqno(2.4)$$

$$ ||f||_{L^\alpha} \le  \sup_{x\ge 0}  \left [K^2_{\alpha/2}(x) x^{1-\nu}\right] \   ||f||_{\nu,1},\eqno(2.5)$$
where $\Gamma(z),\ B(z,w)$ are  Euler's gamma and beta functions, respectively.}
\begin{proof}      In fact,  with the definition of the norm (2.1) and the H$\ddot {o}$lder inequality we obtain
$$||f||_{L^\alpha} = \int_0^\infty K^2_{\alpha/2}(x)|f(x)|dx \le \left(\int_0^\infty K_{\alpha/2} ^{2q} (x)
x^{(1-\nu)q-1}dx\right)^{1/q} ||f||_{\nu, p},\   q={p\over  p-1}\eqno(2.6)$$
and the latter integral via asymptotic behavior of the Macdonald function (1.4), (1.5), (1.6) converges for 
$\nu+\alpha < 1$.   Hence integral (1.7) and  the generalized Minkowski  inequality yield
$$ \left(\int_0^\infty K_{\alpha/2} ^{2q} (x)
x^{(1-\nu)q-1}dx\right)^{1/q} = \left(\int_0^\infty x^{(1-\nu)q-1}
\left(\int_0^\infty e^{- x\cosh u}\cosh (\alpha u /2) du\right)^{2q} dx\right)^{1/q}$$
$$\le \left(\int_0^\infty  \cosh (\alpha u /2) du \left(\int_0^\infty x^{(1-\nu)q-1} e^{- 2q x\cosh u}dx\right)
^{1/q}\right)^2 $$$$ = (2 q)^{2(\nu-1)}\Gamma^{2/q}\left(q(1-\nu)\right)\left( \int_0^\infty { \cosh (\alpha u /2) \over
\cosh^{1-\nu} u} du \right)^2.$$
Calculating the integral with hyperbolic functions via relation (2.4.4.4) in \cite{prud},  Vol. I, we come up with the estimate (2.4).   For the case $p=1$ we end up  immediately  with (2.5), using (2.6), where the supremum is finite via the condition $\nu+\alpha <1$.  Thus the embedding (2.3) is established and Lemma 1 is proved.

\end{proof}

Letting in (1.8) $z=x \in \mathbb{R}_+, \  \mu=i\tau$ and making a  simple substitution,  equality (1.8) becomes
$$\left| K_{(i\tau+ \alpha)/2}(x)\right|^2  = \int_{-\infty} ^\infty K_{i\tau}\left(2x\cosh t \right) e^{\alpha t} dt, \   x >0.\eqno(2.7)$$
Moreover,  appealing to the representation (1.7) of the Macdonald function, we substitute the corresponding integral into the right-hand side of (2.7). Changing the order of integration due to the absolute convergence of the iterated integrals, we find the formula 
$$\Phi_{\alpha,\tau}(x)\equiv \left| K_{(i\tau+ \alpha)/2}(x)\right|^2  =  \int_{-\infty} ^\infty K_{\alpha}\left(2x\cosh t \right)  e^{i\tau t} dt, \   x >0.\eqno(2.8)$$
The representation (2.8) is a key identity, which will be involved to establish a differential  difference equation for the kernel  $\Phi_{\alpha, \tau}(x)$.  Precisely, it has

{\bf Lemma 2}.  {\it The kernel $\Phi_{\alpha, \tau} (x)$ satisfies the following second order differential - difference equation}
$$\frac {d^2 \Phi_{\alpha, \tau}} {dx^2}   + {1\over x}\  {d \Phi_{\alpha, \tau} \over dx} +  {\tau^2\over x^2} \Phi_{\alpha, \tau}= \Phi_{2+ \alpha, \tau} + 2 \Phi_{\alpha, \tau} +  \Phi_{2- \alpha, \tau},\quad x > 0.\eqno(2.9)$$

\begin{proof}  In fact, differentiating two times both sides of (2.8) with respect to $x$, motivating it by the absolute and uniform convergence of the integral and its derivatives, we obtain ($\prime$ means the derivative)
$$ \Phi^{\prime\prime} _{\alpha, \tau} (x)= 4 \int_{-\infty} ^\infty K^{\prime\prime} _{\alpha}\left(2x\cosh t \right)  e^{i\tau t} \cosh^2 t dt=  4 \int_{-\infty} ^\infty K^{\prime\prime} _{\alpha}\left(2x\cosh t \right)  e^{i\tau t} dt$$$$+
4 \int_{-\infty} ^\infty K^{\prime\prime} _{\alpha}\left(2x\cosh t \right)  e^{i\tau t} \sinh^2 t dt.$$  
Meanwhile, the second derivative of the Macdonald function can be expressed as (cf. \cite{erd}, Vol. II)
$$K^{\prime\prime} _{\alpha} (z)= {1\over 4} \left[ K_{2+\alpha}(z) + K_{2-\alpha}(z) \right] + {1\over 2} K_\alpha(z).$$
Hence,
$$ \frac {d^2 \Phi_{\alpha, \tau}} {dx^2}  -  \Phi _{2+\alpha, \tau} (x) -  \Phi _{2-\alpha, \tau} (x)-  2  \Phi _{\alpha, \tau} (x)
= 4 \int_{-\infty} ^\infty K^{\prime\prime} _{\alpha}\left(2x\cosh t \right)  e^{i\tau t} \sinh^2 t dt.\eqno(2.10)$$
On the other hand, integrating twice by parts in the latter integral and eliminating the integrated terms via asymptotic formula (1.4), we find $(\tau \neq 0) $ 
$$4 \int_{-\infty} ^\infty K^{\prime\prime} _{\alpha}\left(2x\cosh t \right)  e^{i\tau t} \sinh^2 t dt = - {2i\tau \over x}  \int_{-\infty} ^\infty K^{\prime} _{\alpha}\left(2x\cosh t \right)  e^{i\tau t} \sinh t dt$$$$-   {2 \over x}  \int_{-\infty} ^\infty K^{\prime} _{\alpha}\left(2x\cosh t \right)  e^{i\tau t} \cosh t dt=  - \left({\tau \over x}\right)^2 \Phi _{\alpha, \tau} (x)
+   {2 \over ix \tau}  \int_{-\infty} ^\infty K^{\prime} _{\alpha}\left(2x\cosh t \right)  e^{i\tau t} \sinh t  dt$$$$
+    {4 \over i \tau}  \int_{-\infty} ^\infty K^{\prime\prime} _{\alpha}\left(2x\cosh t \right)  e^{i\tau t} \cosh t \sinh t  dt
=  - \left({\tau \over x}\right)^2 \Phi _{\alpha, \tau} (x)   - {1 \over x^2} \Phi _{\alpha, \tau} (x)$$$$+ 
{2 \over i \tau} {d\over dx}  \int_{-\infty} ^\infty K^{\prime} _{\alpha}\left(2x\cosh t \right)  e^{i\tau t} \sinh t  dt
=  - \left({\tau \over x}\right)^2 \Phi _{\alpha, \tau} (x)   - {1 \over x^2} \Phi _{\alpha, \tau} (x) - {d\over dx} \left({1 \over x} \Phi _{\alpha, \tau} (x)\right)$$$$=   - \left({\tau \over x}\right)^2 \Phi _{\alpha, \tau} (x)   - {1 \over x} {d \Phi _{\alpha, \tau} (x)\over dx}.$$

Hence, the right-hand side of the latter equality substituting in (2.10), we end up with (2.9).

\end{proof}

{\bf Theorem 1.}   {\it The index transform  $(1.1)$  is well-defined as a  bounded operator from $L^\alpha,   \alpha \ge 0$ into the space $C_0(\mathbb{R})$ of bounded continuous functions vanishing at infinity.   Besides, the following composition representation holds
$$F_\alpha(\tau)=  \left( \mathcal{F} (\mathcal{K}_\alpha f)(2 \cosh t)\right) (\tau),\eqno(2.11)$$
where 
$$ (\mathcal{F} f) (\tau)=   \int_{-\infty} ^\infty f(t)  e^{i\tau t} dt\eqno(2.12)$$
is the operator of Fourier transform and 
$$(\mathcal{K}_\alpha f)(x) = \int_{0} ^\infty K_{\alpha}(xt) f(t) dt\eqno(2.13)$$
is the operator of the Meijer K- transform (cf. \cite{yal})}.

\begin{proof}  In fact, since (see (1.7))  $|K_{(i\tau+\alpha)/2}(x) | \le K_{\alpha/2}(x)$  we have 
$$|F_\alpha(\tau)|  \le  \int_0^\infty K^2_{\alpha/2}(x) | f(x)| dx = ||f||_{L^\alpha}< \infty ,$$
which means that the operator (1.1) is well-defined and the integral converges absolutely and uniformly with respect to $\tau \in \mathbb{R}$.  Thus $F_\alpha (\tau)$ is continuous.   On the other hand, recalling (2.8), we derive
$$|F_\alpha(\tau)|  \le  \int_0^\infty  \int_{-\infty} ^\infty K_{\alpha}\left(2x\cosh t \right)  |f(x)| dx dt =  ||f||_{L^\alpha}< \infty .$$
Hence in view of Fubini's theorem one  can invert the order of integration in the corresponding iterated integral and arrive at the composition (2.11).   Moreover,  the previous estimate says that $(\mathcal{K}_\alpha f)(2\cosh t) \in L_1(\mathbb{R})$. Consequently, $F_\alpha(\tau)$ vanishes at infinity owing to the Riemann-Lebesgue lemma. 
\end{proof} 

{\bf Corollary 1.}  {\it The operator $F_\alpha : L_{\nu,p}(\mathbb{ R}_+) \to L_p(\mathbb{ R}), \ p \ge 2,  \   \alpha+ \nu < 1$ is bounded and }
$$||F_\alpha ||_{L_p(\mathbb{R})} \le \pi^{{1\over p}} 2^{ - 2/q-\nu} q^{\nu-1} \left[\Gamma\left({q\over 2}(1-\nu)\right)\right]^{2/q} $$$$\times B\left({1-\nu+\alpha \over 2},\  {1-\nu-\alpha \over 2} \right)  ||f||_{\nu,p}, \quad  q={p\over p-1}.\eqno(2.14)$$

\begin{proof}  Indeed,  taking the composition (2.11) via Lemma 1 and appealing to the Hausdorff-Young inequality for  Fourier transform (2.12) (cf. \cite{tit},  Theorem 74)
$$|| \mathcal{F} f ||_{L_q(\mathbb{R})} \le (2\pi)^{1/q}  || f ||_{L_p(\mathbb{R})},\  1 < p \le 2,
\ q= {p\over p-1},\eqno(2.15)$$
we find
$$||F_\alpha ||_{L_p(\mathbb{R})}  \le 2 \pi^{{1\over p}}  \left(\int_{0} ^\infty 
\left|(\mathcal{K}_\alpha f)(2 \cosh t)\right|^q dt \right)^{1/q}.\eqno(2.16)$$
Hence by the generalized Minkowski and H$\ddot{o}$lder inequalities  with relation  (2.16.2.2) from \cite{prud},  Vol. II  we obtain similar to the proof of Lemma 1 
$$2 \pi^{{1\over p}}  \left(\int_{0} ^\infty \left|(\mathcal{K}_\alpha f)(2 \cosh t)\right|^q dt \right)^{1/q}
\le 2\pi^{{1\over p}}  \int_0^\infty   |f(x)|\left(\int_0^\infty  K^q_\alpha(2x \cosh t) dt\right)^{1/q}dx$$
$$\le 2 \pi^{{1\over p}}  \int_0^\infty   \int_0^\infty  |f(x)| \cosh(\alpha u) \left(\int_0^\infty  e^{-2qx \cosh t \cosh u} dt\right)^{1/q} du dx$$
$$= 2 \pi^{{1\over p}}  \int_0^\infty   \int_0^\infty  |f(x)| \cosh(\alpha u) \  K_0^{1/q}(2qx\cosh u) du dx$$
$$\le 2 \pi^{{1\over p}}  ||f||_{\nu,p}  \int_0^\infty   \cosh(\alpha u) \left( \int_0^\infty   x^{(1-\nu)q-1} K_0 (2qx\cosh u) 
dx \right)^{1/q} du $$
$$= \pi^{{1\over p}} 2^{1 - 2/q} q^{\nu-1} \left[\Gamma\left({q\over 2}(1-\nu)\right)\right]^{2/q} 
 ||f||_{\nu,p}  \int_0^\infty   \frac{ \cosh(\alpha u)}{\cosh^{1-\nu} u } du $$
$$= \pi^{{1\over p}} 2^{ - 2/q-\nu} q^{\nu-1}\left[\Gamma\left({q\over 2}(1-\nu)\right)\right]^{2/q} 
 B\left({1-\nu+\alpha\over 2},\   {1-\nu-\alpha\over 2}\right)  ||f||_{\nu,p}.$$
Consequently,  combining with (2.16),  we get (2.14).
\end{proof}

The next result tells when operator (1.1) is compact.  

{\bf Theorem 2.} {\it The operator  $F_\alpha  :  L_{\nu,p}(\mathbb{ R}_+) \to L_q(\mathbb{R}), \  1< p \le 2,  \alpha \ge 0,\ \alpha + \nu < 1, \  q=p/(p-1)$ is compact}.

\begin{proof} The proof is based on approximation of the operator (1.1) by a sequence of compact operators of a finite rank with continuous kernels of compact support.  But to achieve this goal, it is sufficient to verify the following Hilbert-Schmidt-type condition
$$ \int_0^\infty\int_{-\infty}^\infty \left |K_{(i\tau+\alpha)/2}(x) \right|^{2q}
x^{(1-\nu)q-1}d\tau dx < \infty.\eqno(2.17)$$
Indeed, recalling again integral representation (1.7),  (2.8),  the Hausdorff-Young inequality (2.14) and   the generalized Minkowski inequality , we deduce
$$\left( \int_0^\infty\int_{-\infty}^\infty \left |K_{(i\tau+\alpha)/2}(x) \right|^{2q} 
x^{(1-\nu)q-1}d\tau dx \right)^{1/q} $$$$\le  2 \pi^{1/q} \left( \int_0^\infty x^{(1-\nu)q-1} \left( \int_{0}^\infty 
K^p _{\alpha}(2x\cosh t) dt \right)^{1/ (p-1)} dx   \right)^{1/q} $$
$$\le  2 \pi^{1/q} \left( \int_0^\infty x^{(1-\nu)q-1} \left( \int_{0}^\infty  \cosh(\alpha u) 
\  K^{1/p} _{0}(2xp \cosh u)\  du\right)^q  dx   \right)^{1/q} $$
$$\le  2 \pi^{1/q} \int_0^\infty  \cosh(\alpha u) \left( \int_0^\infty x^{(1-\nu)q-1} 
\  K^{1/(p-1)} _{0}(2xp \cosh u)\  dx\right)^{1/q}  du $$
$$= 2^{\nu} \pi^{1/q}  p^{\nu-1} \left( \int_0^\infty x^{(1-\nu)q-1} 
\  K^{1/(p-1)} _{0}(x)\  dx\right)^{1/q}  \int_0^\infty  \frac{ \cosh(\alpha u)}{\cosh^{1-\nu} u } du $$

$$= { \pi ^{1/q}  p^{\nu-1}\over 2}   B\left({1-\nu+\alpha\over 2},\   {1-\nu-\alpha\over 2}\right)\left( \int_0^\infty x^{(1-\nu)q-1} \  K^{1/(p-1)} _{0}(x)\  dx\right)^{1/q}  $$

$$ \le { \pi ^{1/q}  p^{\nu-1}\over  2}   B\left({1-\nu+\alpha\over 2},\   {1-\nu-\alpha\over 2}\right)$$$$\times \left( \int_0^\infty du \left(\int_0^\infty x^{(1-\nu)q-1}  e^{-x\cosh u /(p-1)} dx\right)^{p-1} \right)^{1/p}  $$

$$=  { \pi ^{1/q}  q^{\nu-1} \Gamma^{1/q} (q(1-\nu)) \over 2 \   \Gamma(1-\nu)}   \Gamma\left({1-\nu+\alpha\over 2}\right) \Gamma\left({1-\nu-\alpha\over 2}\right) \left(\int_0^\infty {du\over \cosh^{p(1-\nu)} u }\right)^{1/p}  $$

$$ =  { \pi ^{1/q}  (q/2)^{\nu-1} 2^{-1 -2/p}\   \Gamma^{1/q} (q(1-\nu)) \over   \Gamma^{1/p} (p(1-\nu)) }   \Gamma^{2/p}\left({p\over 2}(1-\nu)\right)  $$$$\times B\left({1-\nu+\alpha\over 2}, \  {1-\nu-\alpha\over 2}\right) < \infty.$$

\end{proof} 

Letting $p=q=2, \alpha =0, \nu= 1/2$ we get 

{\bf Corollary  2.} {\it The Lebedev operator $(1.1)$  $F_0  :  L_{2}(\mathbb{ R}_+) \to L_2(\mathbb{R})$ is the Hilbert - Schmidt operator with  the square of the Macdonald function $K^2_{i\tau/2} (x)$ as the Hilbert-Schmidt kernel. Moreover, its norm is equal to  $\pi^2/2$.}

\begin{proof} Employing the Parseval equality for the Fourier transform \cite{tit}, representation (2.8)  and relation (2.16.33.2) in \cite{prud}, Vol. II, we derive (see (2.17)) 
$$||F_0||=  \left( \int_0^\infty\int_{-\infty}^\infty K^4_{i\tau/2}(x) d\tau dx \right)^{1/2} = 2 \sqrt{\pi} \left( \int_0^\infty
\int_{0}^\infty K^2_{0}(2x\cosh t) dt dx \right)^{1/2}$$
$$ = \sqrt{2 \pi} \left( \int_0^\infty {dt\over \cosh t}  \int_{0}^\infty K^2_{0}(x) \  dx \right)^{1/2} = {\pi^2\over 2}.$$

\end{proof}

Another representation of the transform (1.1) can be given via the Parseval equality for the Mellin transform \cite{tit}
$$\int_0^\infty f(x) g(x) dx= {1\over 2\pi i} \int_{\nu- i\infty}^{\nu+i\infty} f^*(s) g^*(1-s) ds,\eqno(2.18)$$
where $f \in L_{\nu, p}(\mathbb{R}_+),\ 1 < p \le 2$ and 
$$f^*(s)= \int_0^\infty f(x) x^{s-1} dx\eqno(2.19)$$
is its Mellin transform and integral (2.19) converges in mean with respect to the norm in $L_q(\nu- i\infty, \nu + i\infty),\ 
 q=p/(p-1)$.   The inverse Mellin transform is given accordingly
 $$f(x)= {1\over 2\pi i}  \int_{\nu- i\infty}^{\nu+i\infty} f^*(s)  x^{-s} ds,\eqno(2.20)$$
where the integral converges in mean with respect to the norm (2.2) in   $L_{\nu, p}(\mathbb{R}_+)$.    An immediate consequence of Theorems 86, 87 in \cite{tit} is the following result.

{\bf Theorem 3}. {\it Let $f \in L_{\nu, p}(\mathbb{R}_+),\ 1 < p \le 2,  \alpha+\nu < 1.$  Then for all $\tau \in \mathbb{R}$}
$$F_\alpha(\tau) =  {1\over 8i\sqrt \pi }  \int_{\nu- i\infty}^{\nu+i\infty} \Gamma\left({1- s+i\tau\over 2}\right) 
\Gamma\left({1- s-i\tau\over 2}\right)  $$$$\times \frac{\Gamma ((1- s+\alpha)/2) \Gamma ((1- s-\alpha)/2)}
{\Gamma ((1-s)/2) \Gamma (1-   s/ 2)} f^*(s)  ds.\eqno(2.21)$$

\begin{proof} In fact, the proof is based on the equality (2.18) and relation (8.4.23.31) in \cite{prud}, Vol.  III, which drives us to the following representation of the kernel  $\left| K_{(i\tau+ \alpha)/2}(x)\right|^2$
$$\left| K_{(i\tau+ \alpha)/2}(x)\right|^2 =   {1\over 8i\sqrt \pi }  \int_{\mu- i\infty}^{\mu+i\infty} \Gamma\left({s+i\tau\over 2}\right) \Gamma\left({s-i\tau\over 2}\right) \frac{\Gamma ((s+\alpha)/2) \Gamma ((s-\alpha)/2)}
{\Gamma (s/2) \Gamma ((s+1)/2)} x^{-s} ds, \ x >0,\eqno(2.22)$$
where $\tau \in \mathbb{R}$ and $\mu > \alpha$. 
\end{proof}

Finally in this section we investigate the existence and boundedness of the adjoint operator (1.2). In fact, following the general operator theory it can be  established from the boundedness of the operator (1.1). However, we will prove it directly, getting an explicit estimation of its norm.  Assuming $g(\tau) \in L_p(\mathbb{R}),\ 1 <p \le 2$ and recalling (2.8) with the asymptotic formula (1.4) for the Macdonald function, we find that for each $x >0$ the function $K_\alpha(2x\cosh t) \in  L_p(\mathbb{R}),\ 1 <p \le 2$.    Hence via the Parseval theorem for the Fourier transform (cf. \cite{tit}, Theorem 75), operator (1.2) can be written as
$$G_\alpha(x) = \int_{-\infty}^\infty K_\alpha(2x\cosh t)  \left(\mathcal{F} g\right) (t)dt,\   x >0,\eqno(2.23)$$
where $\left(\mathcal{F} g\right) (t) \in L_q(\mathbb{R}),\  q= {p\over p-1}$ is the Fourier transform (2.12) of $g$.

{\bf Theorem 4.}   {\it Let $g \in L_p(\mathbb{R}),\   1 < p \le 2,\    0\le \alpha < 1+ {1\over 2p}$. Then operator  $(1.2)$  is well-defined  and for all $x >0$}
$$|G_\alpha(x)|  \le  \pi^{(1+ 1/q)/2 }  2^{-2 - 1/p} p^{-1/(2p)}   B\left({1+\alpha\over 2}+ {1\over 4p},\   {1-\alpha\over 2}- {1\over 4p} \right)   x^{-1 -1/(2p)}  || g ||_{L_p(\mathbb{R})}.\eqno(2.24)$$

\begin{proof}   Taking  (2.23), we recall the H$\ddot{o}$lder inequality,  the Hausdorff-Young inequality (2.15) and the generalized Minkowski inequality to obtain 
$$|G_\alpha(x)| \le \left(\int_{-\infty}^\infty K^p_\alpha(2x\cosh t) dt \right)^{1/p}  || \mathcal{F} g ||_{L_q(\mathbb{R})} 
\le   \pi^{1/q}    || g||_{L_p(\mathbb{R})} \int_{0}^\infty \cosh(\alpha u) K^{1/p}_0 (2x p \cosh u) du   $$
$$\le  \pi^{1/q}    || g ||_{L_p(\mathbb{R})} \int_{0}^\infty \cosh(\alpha u)\    e^{-2x\cosh u} \left(
\int_{0}^\infty  e^{-2x p \cosh u \  t^2} dt \right)^{1/p}  du $$
 $$=  \pi^{(1+ 1/q)/2 }  2^{- 3/(2p)} (xp)^{-1/(2p)}   || g ||_{L_p(\mathbb{R})} \int_{0}^\infty {\cosh(\alpha u) 
 \over \cosh^{1/(2p)} u } e^{-2x\cosh u}   du $$
$$\le \pi^{(1+ 1/q)/2 }  2^{-1 - 3/(2p)} p^{-1/(2p)}  x^{-1 -1/(2p)}   || g||_{L_p(\mathbb{R})} \int_{0}^\infty {\cosh(\alpha u) 
 \over \cosh^{1+ 1/(2p)} u }  du $$
 $$=  \pi^{(1+ 1/q)/2 }  2^{-2 - 1/p} p^{-1/(2p)}  x^{-1 -1/(2p)} B\left({1+\alpha\over 2}+ {1\over 4p},\   {1-\alpha\over 2}- {1\over 4p} \right) || g ||_{L_p(\mathbb{R})} , $$
which proves (2.24). 

\end{proof} 

{\bf Theorem  5.}  {\it The operator $G_\alpha : L_p(\mathbb{ R}) \to L_{\nu,r}(\mathbb{ R}_+), \   1< p \le 2,  \    r \ge 1,  \   \alpha <  \nu $ is bounded and }
$$||G_\alpha ||_{\nu,r} \le \pi^{1- 1/p} 2^{\nu- 2- 2/p}  {\Gamma^{1/r} (\nu r)\over r^\nu \Gamma^{1/p} (\nu p)}  \  \Gamma^{2/p} \left( {\nu p\over 2}\right)   B \left( {\nu +\alpha\over 2},\   {\nu -\alpha \over 2}\right)  \    || g ||_{L_p(\mathbb{R})}. .$$

\begin{proof}  Indeed,   recalling (2.23),  we apply again the generalized Minkowski,   H$\ddot{o}$lder inequalities and  the Hausdorff-Young inequality   (2.15) to find
$$||G_\alpha ||_{\nu,r}  \le \int_{-\infty}^\infty \left|  \left(\mathcal{F} g\right) (t)\right|  \left(\int_{0} ^\infty 
x^{\nu r -1} K^r_\alpha(2x\cosh t) dx \right)^{1/r}  dt$$$$\le  || \mathcal{F} g||_{L_q(\mathbb{R})}
\left( \int_{-\infty}^\infty  \left(\int_{0} ^\infty  x^{\nu r -1} K^r_\alpha(2x\cosh t) dx \right)^{p/r}  dt\right)^{1/p}$$
$$\le (2\pi)^{1/q} 2^{-\nu}  || g ||_{L_p(\mathbb{R})} \left( \int_{-\infty}^\infty  {dt\over \cosh^{\nu p} t} \right)^{1/p} 
\left(\int_{0} ^\infty  x^{\nu r -1} K^r_\alpha(x) dx \right)^{1/r} $$
$$= \pi^{1/q} 2^{-2/p}  B^{1/p} \left( {\nu p\over 2},\   {\nu p\over 2}\right)  
\left(\int_{0} ^\infty  x^{\nu r -1} K^r_\alpha(x) dx \right)^{1/r}  || g ||_{L_p(\mathbb{R})}  $$
$$\le \pi^{1/q} 2^{-2/p}  B^{1/p} \left( {\nu p\over 2},\   {\nu p\over 2}\right)  
\int_0^\infty  \cosh(\alpha t) \left(\int_{0} ^\infty  x^{\nu r -1} e^{-xr \cosh u} dx  \right)^{1/r} dt  \   || g ||_{L_p(\mathbb{R})}  $$
$$= \pi^{1/q} 2^{-2/p} r^{-\nu} \Gamma^{1/r} (\nu r)  B^{1/p} \left( {\nu p\over 2},\   {\nu p\over 2}\right)  
\    || g ||_{L_p(\mathbb{R})} \int_0^\infty  {\cosh(\alpha t) \over \cosh^\nu t } dt   $$
$$= \pi^{1/q} 2^{\nu- 2- 2/p}  {\Gamma^{1/r} (\nu r)\over r^\nu \Gamma^{1/p} (\nu p)}  \  \Gamma^{2/p} \left( {\nu p\over 2}\right)   B \left( {\nu +\alpha\over 2},\   {\nu -\alpha \over 2}\right)  \    || g ||_{L_p(\mathbb{R})}. $$ 

\end{proof}

\section{Inversion theorems}   

The composition representation (2.11) and the properties of the Fourier and Mellin transforms are key ingredients to prove the inversion theorem for the index transform (1.1).   Namely, we have

{\bf Theorem 6}.  {\it  Let $0<  \alpha <  1, \     \nu < 1-\alpha,\  1< p \le 2,\  q= p/(p-1). $    Let $s f^*(s) \in L_p(\nu-i\infty, \nu+i\infty)$, where $f^*(s)$ is the Mellin transform $(2.19)$ of $f \in L_1((1,\infty); \  t^\alpha dt)$, i.e. $f(t)$ is integrable over $(1,\infty)$ with respect to the measure $t^\alpha  dt$. If, besides above assumptions,  the generalized Lebedev  transform $(1.1)$ of $f$ satisfies the condition $\tau e^{\pi |\tau|}  F_\alpha(\tau) \in L_1(\mathbb{R})$ and its Mellin transform vanishes at the point $1-\alpha$, i.e. $f^*(1-\alpha)=0$, then for all $x> 0$ the following inversion formula holds
$$f(x) =  { 1\over  \pi  }   \int_{-\infty} ^\infty \left[  \frac{   (x/2)^{i\tau-1}}{\Gamma\left((i\tau- \alpha)/ 2\right)  \Gamma\left((\alpha+ i\tau)/ 2\right)}  \  {}_2F_3\left( {i\tau\over 2},\    {1+ i\tau\over 2};  \ 1+ i\tau,   {i\tau-\alpha \over 2}, \   {\alpha+ i\tau\over  2} ;  \   x^2 \right)\right. $$
$$\left.   -  \   \frac{2\cosh(\pi\tau/2) }{x \  \Gamma(\alpha/2)\Gamma(-\alpha/2)}\right] \      F_\alpha(\tau)\  d\tau, \eqno(3.1)$$
where ${}_2F_{3}(a_1, \  a_2; \  b_1, \  b_2,\  b_3;  z )$ is the generalized hypergeometric function and the integral converges absolutely.}

\begin{proof}  In fact, since $s f^*(s) \in L_p(\nu-i\infty, \nu+i\infty),\ 1 < p\le 2$,  it means that  $f^*(s) \in L_p(\nu-i\infty, \nu+i\infty)$. Hence Theorem 86 in \cite{tit} says that $f$, which is  given by formula (2.20),   belongs to $L_{\nu,q}(\mathbb{R}_+)$.  Then by virtue of Lemma 1 and Theorem 1 we observe that $F_\alpha(\tau)$ is continuous.  Therefore the condition $\tau e^{\pi |\tau|}  F_\alpha(\tau) \in L_1(\mathbb{R})$  implies  $F_\alpha \in L_1(\mathbb{R})$.  Hence (2.11) and the inverse Fourier transform yield the equality 
$$ {1\over 2\pi} \int_{-\infty} ^\infty F_\alpha(\tau)  e^{i\tau t} d\tau =  \int_{0} ^\infty K_{\alpha}(2x\cosh t) f(x) dx,$$
and after simple substitution $\lambda =\cosh t$ it becomes 
$$ {1\over 2\pi} \int_{-\infty} ^\infty F_\alpha(\tau)  e^{i\tau \log\left(\lambda + \sqrt{\lambda^2-1} \right)} d\tau =  \int_{0}^\infty K_{\alpha}(2x\lambda ) f(x) dx, \  \lambda  \ge 1.\eqno(3.2)$$
The integral in the left-hand side of (3.2) converges absolutely.   Further,  recalling the H$\ddot{o}$lder inequality, asymptotic formulas (1.4), (1.5) for the Macdonald function and the condition $f \in L_{\nu,q}(\mathbb{R}_+)$, it is not difficult to verify the absolute convergence of the integral in the right-hand side of (3.2).    Moreover, it permits a differentiation with respect to $\lambda$ in (3.2). As a result we obtain
$$ {i \over 2\pi} \int_{-\infty} ^\infty \tau F_\alpha(\tau)  {e^{i\tau \log\left(\lambda + \sqrt{\lambda^2-1} \right)}\over 
\sqrt{\lambda^2-1} }  d\tau =  2 \int_{0} ^\infty K^\prime _{\alpha}(2x\lambda )  x f(x) dx.\eqno(3.3)$$
Integrating by parts in the right-hand side of (3.3) and eliminating the integrated terms, we get
$${i \over 2\pi} \int_{-\infty} ^\infty \tau F_\alpha(\tau)  {e^{i\tau \log\left(\lambda + \sqrt{\lambda^2-1} \right)}\over 
\sqrt{\lambda^2-1} }  d\tau =  -  {1\over \lambda}  \int_{0} ^\infty K_{\alpha}(2x\lambda )  {d\over dx} \left[x f(x)\right] dx.\eqno(3.4)$$
But since $s f^*(s) \in L_p(\nu-i\infty, \nu+i\infty)$,  one has that $f(x)$ is equivalent to some absolutely continuous function,  ${d\over dx} \left[x f(x)\right]  \in L_{\nu,q}(\mathbb{R}_+)$ and 
$${d\over dx} \left[x f(x)\right] =  {1\over 2\pi i}  \int_{\nu- i\infty}^{\nu+i\infty} (1-s) f^*(s)  x^{-s} ds,$$
where the integral converges in mean with respect to the norm  in   $L_{\nu, q}(\mathbb{R}_+)$.  Further, the Parseval equality (2.18) and relation (2.16.2.2) in \cite{prud}, Vol. II allow us to write for all $\lambda \ge 1$ 
$$ \int_{0} ^\infty K_{\alpha}(2x\lambda )  {d\over dx} \left[x f(x)\right] dx =  {1\over 8\pi i}  \int_{\nu- i\infty}^{\nu+i\infty} (1-s) \  \Gamma\left({1-s+\alpha\over 2}\right)  \Gamma\left({1-s-\alpha\over 2}\right)   f^*(s) \ \lambda^{s-1} ds.$$
Therefore, combining with (3.4), we have
$${i \over 2\pi} \int_{-\infty} ^\infty \tau F_\alpha(\tau)  {e^{i\tau \log\left(\lambda + \sqrt{\lambda^2-1} \right)}\over 
\sqrt{\lambda^2-1} }  d\tau = -   {1\over 8\pi i}  \int_{\nu- i\infty}^{\nu+i\infty} (1-s) \  \Gamma\left({1-s+\alpha\over 2}\right)  \Gamma\left({1-s-\alpha\over 2}\right) $$$$\times   f^*(s) \ \lambda^{s-2} ds.\eqno(3.5)$$
Meanwhile, $F_\alpha(\tau)$ is even.   Thus it implies from (3.5)
$${1\over \pi} \int_{-\infty} ^\infty \tau F_\alpha(\tau)  {\sin\left(\tau \log\left(\lambda + \sqrt{\lambda^2-1} \right)\right)\over 
\sqrt{\lambda^2-1} }  d\tau =   {1\over 4\pi i}  \int_{\nu- i\infty}^{\nu+i\infty} (1-s) \  \Gamma\left({1-s+\alpha\over 2}\right)  \Gamma\left({1-s-\alpha\over 2}\right) $$$$\times   f^*(s) \ \lambda^{s-2} ds.\eqno(3.6)$$
In the meantime, relations (2.16.2.2),  (2.16.6.1) in \cite{prud},  Vol. II and the Parseval equality (2.18) give the following representations of the kernel in the left-hand side of (3.6)
$$ {\sin\left(\tau \log\left(\lambda + \sqrt{\lambda^2-1} \right)\right)\over \sqrt{\lambda^2-1} } = {1\over \pi} \sinh(\pi\tau) \int_0^\infty e^{-\lambda x} K_{i\tau}(x) dx$$$$=  {\sinh(\pi\tau) \over \pi^2 i}    \int_{\mu- i\infty}^{\mu+i\infty} 2^{-2-s} \Gamma(s)  \  \Gamma\left({1-s+ i\tau\over 2}\right)  \Gamma\left({1-s- i\tau\over 2}\right) \lambda^{-s} ds,\ 0 < \mu < 1.\eqno(3.7)$$
Hence, substituting the right-hand side of  (3.7) into the left-hand side of (3.6) and changing the order of integration via Fubini's theorem due to the estimate
$$ \int_{-\infty} ^\infty \left| \tau F_\alpha(\tau) \sinh(\pi\tau) \right|  \int_{\mu- i\infty}^{\mu+i\infty} \left|2^{-2-s} \Gamma(s)  \  \Gamma\left({1-s+ i\tau\over 2}\right)  \Gamma\left({1-s- i\tau\over 2}\right) \lambda^{-s} ds\right|   d\tau $$
$$\le {\Gamma^2 ((1-\mu)/2) \over 2^{2+\mu} \lambda^\nu} \int_{-\infty} ^\infty \left|\tau  F_\alpha(\tau) \right| e^{\pi|\tau|} d\tau   \int_{\mu- i\infty}^{\mu+i\infty} \left| \Gamma(s)  ds \right|< \infty, \eqno(3.8)$$ 
it becomes by virtue of simple changes of variables
$${1\over \pi^3 i}  \int_{\mu- i\infty}^{\mu+i\infty} 2^{-s-2} \Gamma(s)   \int_{-\infty} ^\infty \tau \sinh(\pi\tau) \   F_\alpha(\tau) \Gamma\left({1-s+ i\tau\over 2}\right)  \Gamma\left({1-s- i\tau\over 2}\right)  d\tau  \    \lambda^{-s} ds $$$$=   {1\over 4\pi i}  \int_{2-\nu- i\infty}^{2-\nu+i\infty} (s-1) \  \Gamma\left({s-1+\alpha\over 2}\right)  \Gamma\left({s-1 -\alpha\over 2}\right)   f^*(2-s) \ \lambda^{-s} ds.\eqno(3.9)$$
In the meantime, we will show that under conditions of the theorem $f^*(s)$ is analytic in the strip $\nu  < {\rm Re} \ s < 1+\alpha$. Indeed, it follows immediately from the absolute and uniform convergence of the integral (2.19).   We have $(s=\mu+i\tau)$,
$$|f^*(s)| \le \int_0^1 | f(t)| t^{\mu-1} dt +  \int_1^\infty  | f(t)| t^{\mu-1} dt \le \left( \int_0^1 | f(t)|^q  t^{\nu q-1} dt\right)^{1/q}
 \left( \int_0^1  t^{p(\mu- \nu)-1} dt\right)^{1/p} $$$$+  \int_1^\infty  | f(t)|\  t^\alpha  dt \le [p(\mu-\nu)]^{-1/p}  ||f||_{\nu,q} +
 \int_1^\infty  | f(t)|\  t^\alpha  dt ,\eqno(3.10)$$
and the result follows.  Moreover,  appealing to the Stirling asymptotic formula for  gamma- functions \cite{erd}, we observe that the integrand in the right-hand side of (3.9) belongs to $L_1(\mu-i\infty, \mu +i\infty)$ for all $\mu \in (1-\alpha,\  2-\nu]$.   Therefore, the Cauchy theorem and the condition $f^*(1-\alpha)=0$ permits us to shift the contour $(2-\nu-i\infty, 2-\nu+i\infty)$ to the left, making the  integration  over the line $(\mu-i\infty, \mu +i\infty)$ with $1-\alpha < \mu < 1$.  Hence   we arrive at the equality
$${1\over \pi^3 i}  \int_{\mu - i\infty}^{\mu +i\infty} 2^{-s-2} \Gamma(s)  \    \lambda^{-s}   \int_{-\infty} ^\infty \tau \sinh(\pi\tau)    F_\alpha(\tau) \Gamma\left({1-s+ i\tau\over 2}\right)  \Gamma\left({1-s- i\tau\over 2}\right)  d\tau  \ ds $$$$=   {1\over 4\pi i}  \int_{\mu- i\infty}^{\mu +i\infty} (s-1) \  \Gamma\left({s-1+\alpha\over 2}\right)  \Gamma\left({s-1 -\alpha\over 2}\right)   f^*(2-s) \ \lambda^{-s} ds.\eqno(3.11)$$
Further,   taking into account (3.10) and estimating the integrand in the left-hand side of (3.11) as a function of $s$ similar to (3.8), we observe that both integrands in (3.11) are from the space $L_1(\mu -i\infty,  \mu + i\infty),\ \mu \in (1-\alpha, 1)$. Consequently, the uniqueness theorem for the Mellin transform (see \cite{tit}) immediately drives us at the equality
$${2^{-s} \Gamma(s) \over \pi^2 }   \int_{-\infty} ^\infty \tau \sinh(\pi\tau)  \   F_\alpha(\tau) \Gamma\left({1-s+ i\tau\over 2}\right)  \Gamma\left({1-s- i\tau\over 2}\right)  d\tau  $$$$=  (s-1) \  \Gamma\left({s-1+\alpha\over 2}\right)  \Gamma\left({s-1 -\alpha\over 2}\right)   f^*(2-s),\quad    s \in (\mu -i\infty,  \mu + i\infty).$$
Hence the simple change of variables  and the reduction formula for the gamma-function yield 
$$ f^*(s) = {2^{s-2} \Gamma(1- s) \over \pi^2  }   \int_{-\infty} ^\infty \tau \sinh(\pi\tau)  \frac{ \Gamma\left(( s-1 + i\tau)/ 2\right)  \Gamma\left((s-1 - i\tau)/ 2\right) } {  \Gamma\left((1-s+\alpha)/ 2\right)  \Gamma\left((1-s -\alpha)/ 2\right)} \  F_\alpha(\tau)  d\tau, \eqno(3.12)$$
where $s \in (2-\mu-i\infty,  2- \mu + i\infty)$,  and it is not difficult to verify the condition $f^*(1-\alpha)=0$.   But the right-hand side of (3.12) is integrable over  the line $(2-\mu-i\infty,  2- \mu + i\infty)$.  Indeed,  we have 
$$\int_{2-\mu - i\infty}^{2- \mu +i\infty} \left|  \Gamma(1- s)   \int_{-\infty} ^\infty \tau \sinh(\pi\tau)  \frac{ \Gamma\left(( s-1 + i\tau)/ 2\right)  \Gamma\left((s-1 - i\tau)/ 2\right) } {  \Gamma\left((1-s+\alpha)/ 2\right)  \Gamma\left((1-s -\alpha)/ 2\right)} \  F_\alpha(\tau)  d\tau ds \right| $$
$$= \int_{2-\mu - i\infty}^{2- \mu +i\infty} \left|  \frac{ \Gamma(1- s) \Gamma(s-1) }{ \Gamma\left((1-s+\alpha)/ 2\right)  \Gamma\left((1-s -\alpha)/ 2\right)}  \right.$$$$\left. \times \int_{-\infty} ^\infty \tau \sinh(\pi\tau) B \left(( s-1 + i\tau)/ 2,\  (s-1 - i\tau)/ 2\right)  \  F_\alpha(\tau)  d\tau ds \right|$$
$$\le   B(1-\mu,\ 1-\mu) \int_{2-\mu - i\infty}^{2- \mu +i\infty} \left|  \frac{ \Gamma(1- s) \Gamma(s-1) }{ \Gamma\left((1-s+\alpha)/ 2\right)  \Gamma\left((1-s -\alpha)/ 2\right)}  ds \right| \int_{-\infty} ^\infty |\tau  \  F_\alpha(\tau) | e^{\pi|\tau|}  d\tau < \infty\eqno(3.13)$$
via the condition  $\tau e^{\pi |\tau|}  F_\alpha(\tau) \in L_1(\mathbb{R})$  and the estimate 
$$ \frac{ \Gamma(1- s) \Gamma(s-1) }{ \Gamma\left((1-s+\alpha)/ 2\right)  \Gamma\left((1-s -\alpha)/ 2\right)}  
= O\left(|u|^{1-\mu} e^{-\pi|u|/2}\right),  \quad s= 2-\mu+iu, \   |u| \to \infty.$$
Consequently,  applying to both sides of (3.12) the inverse Mellin transform (2.20), we come up with the equality for all $x >0$
$$ f(x) =  {1\over 2 \pi^2 \sqrt \pi}   \int_{-\infty} ^\infty \tau \sinh(\pi\tau) S_\alpha (x,\tau) F_\alpha(\tau) d\tau, \ x >0,$$
where under a simple substitution with the use of the reduction formula for the gamma function
$$S_\alpha(x,\tau) = {1\over 2 \pi i}  \int_{1- \mu/2 - i\infty}^{1- \mu/2 +i\infty}   \frac{ \Gamma\left(s- (1+ i\tau)/ 2\right)  \Gamma\left(s- (1- i\tau)/ 2\right) \Gamma(3/2 -s) \Gamma(1- s)}{\Gamma\left((1+\alpha)/ 2- s\right)  \Gamma\left((1 -\alpha)/ 2- s\right)(1/2- s)} x^{-2s} ds,\eqno(3.14)$$
and  the change of the corresponding order of integration  is allowed owing to the estimate (3.13) and Fubini's  theorem.
Our goal now is  to calculate the kernel $S_\alpha(x,\tau)$ in terms of the generalized hypergeometric functions ${}_2F_3$ via the Slater theorem (cf. in \cite{prud}, Vol. III).  In fact,  the integral (3.14) is equal to 
$$S_\alpha(x,\tau) =   x^{-(1+i\tau)} \sum_{n=0}^\infty {(-1)^n x^{2n}\over n!}  
\frac{\Gamma(i\tau-n)  \Gamma(n - i\tau/2) \Gamma(n +(1- i\tau) /2)}{\Gamma\left(n +(\alpha- i\tau)/ 2\right)  \Gamma\left(n-  (i\tau +\alpha)/ 2\right)} $$$$+    x^{i\tau-1} \sum_{n=0}^\infty {(-1)^n x^{2n}\over n!}  
\frac{\Gamma(- i\tau-n)  \Gamma(n + i\tau/2) \Gamma(n +(1+ i\tau) /2)}{\Gamma\left(n +(\alpha+ i\tau)/ 2\right)  \Gamma\left(n-  (\alpha-i\tau )/ 2\right)}$$
$$ -  \  \frac{2\pi \sqrt \pi }{x \tau\sinh(\pi\tau/2) \Gamma(\alpha/2)\Gamma(-\alpha/2)}  = {\pi  \  \sqrt\pi   \over  \tau \sinh(\pi\tau)}  \left[  \frac{  \ (x/2)^{i\tau-1}}{\Gamma\left((i\tau- \alpha)/ 2\right)  \Gamma\left((\alpha+ i\tau)/ 2\right)} \right.$$
$$\times 
{}_2F_3\left( {i\tau\over 2},\    {1+ i\tau\over 2};  \ 1+ i\tau,    \   {\alpha+ i\tau\over  2},\    { i\tau - \alpha\over  2} ;  \   x^2 \right) +  \frac{ (2/x)^{1+i\tau} }{\Gamma\left(-(\alpha+i\tau)/ 2\right)  \Gamma\left((\alpha- i\tau)/ 2\right)}$$
$$\left. \times  {}_2F_3\left( - {i\tau\over 2},\    {1-i\tau\over 2};  \ 1-i\tau,   \   {\alpha- i\tau\over  2},\  -  { i\tau +\alpha\over  2} ;  \   x^2 \right)   -  \   \frac{4\cosh(\pi\tau/2) }{x  \Gamma(\alpha/2)\Gamma(-\alpha/2)}\right].$$
Finally,  since $F_\alpha$ is even, we easily arrive  at the inversion formula (3.1), completing the proof of the theorem.

\end{proof} 

In order to prove the  corresponding inversion theorem  for  the index transform (1.2), we will need  the value of the integral, involving   squares of modules of the gamma functions.    Precisely,  we have

{\bf  Lemma 3}. {\it Let,  $\varepsilon >0, \  x \in \mathbb{R}.$ Then}
$$\int_{-\infty}^\infty  \left| \Gamma\left({\varepsilon + i(x- t)\over 2}\right)  \Gamma\left({\varepsilon + i(x+ t)\over 2}\right) \right|^2 dt=  4 \pi \left|\Gamma(\varepsilon+ ix)\right|^2 B(\varepsilon, \varepsilon) .\eqno(3.15)$$

\begin{proof} The proof  is based on the Parseval equality for the Fourier cosine transform \cite{tit} and  relation (2.5.46.6) in \cite{prud}, Vol.  I.   Thus we obtain
$$\int_{-\infty}^\infty  \left| \Gamma\left({\varepsilon + i(x- t)\over 2}\right)  \Gamma\left({\varepsilon + i(x+ t)\over 2}\right) \right|^2 dt =  2^{4-2\varepsilon} \pi \left|\Gamma(\varepsilon+ ix)\right|^2 \int_0^\infty {dt\over \cosh^{2\varepsilon} t }
$$$$= 4 \pi \left|\Gamma(\varepsilon+ ix)\right|^2 B(\varepsilon, \varepsilon).$$

\end{proof}

{\bf Definition 1}.  {\it A function $f$ satisfies the Dini condition at some point $x$, if 
$$\varphi(t)= {f(x+t)- f(x)\over t} $$
is integrable in some neighborhood  of the origin.}

Now the inversion formula for the index transform (1.2) is given by

{\bf Theorem 7.} {\it Let $\alpha \in \mathbb{R},\  g(\tau) \in L_p(\mathbb{R}),\ 1 < p \le 2,$  satisfying the Dini condition at some point $x \in \mathbb{R}\backslash \{0\} $. Let  also the index transform $G_\alpha$ be such that its Mellin transform $ G^*_\alpha(s) \in L_1((\mu-i\infty,\  \mu +i\infty);  |s|^{1/2} ds ),\   \mu >0$.  Then 
$$ g(x) = \lim_{\varepsilon \to 0+} {1\over \pi} \  \frac{ \Gamma(2\varepsilon) }{ \Gamma(\varepsilon)}
\int_0^\infty \left[\frac{  \Gamma(-ix)\left(t/2\right) ^{\varepsilon +ix-1} }{ \Gamma(\varepsilon- ix) \Gamma\left((\varepsilon+\alpha+  ix)/2 \right) \Gamma\left((\varepsilon +  ix-\alpha)/2\right)} \right.$$$$\left. \times  {}_2F_3\left( {\varepsilon + ix\over 2},\    {\varepsilon+  ix + 1\over 2};  \ 1+ ix,    \   {\varepsilon+ \alpha+ ix\over  2},\    { \varepsilon+ ix - \alpha\over  2} ;  \   t^2 \right)   + \frac{  \Gamma(ix)\left(t/ 2\right)^{\varepsilon - ix-1} }{ \Gamma(\varepsilon+ ix) \Gamma\left((\varepsilon+\alpha-  ix)/2 \right) \Gamma\left((\varepsilon -  ix-\alpha)/2\right)} \right.$$$$\left. \times  {}_2F_3\left( {\varepsilon - ix\over 2},\    {\varepsilon-  ix + 1\over 2};  \ 1- ix,    \   {\varepsilon+ \alpha- ix\over  2},\    { \varepsilon- ix - \alpha\over  2} ;  \   t^2 \right) \right]G_\alpha(t) dt, \eqno(3.16)$$
where the limit is pointwise.} 

\begin{proof}   Recalling formula (2.23), which is valid under conditions of the theorem, we calculate the Mellin transform 
(2.19) from its both sides with $s$ belongs to the vertical line with ${\rm Re}\ s > |\alpha|$.  After the change of the order of integration via Fubini's theorem by virtue of the absolute convergence and the use of relation (2.16.2.2) in \cite{prud}, Vol. II, we find 
$$\frac{4\  G_\alpha^*(s)}{ \Gamma\left((s + \alpha)/2\right)  \Gamma\left((s- \alpha)/2\right)}= 
  \int_{-\infty}^\infty  {\left(\mathcal{F} g\right) (t)\over \cosh^s t} dt.\eqno(3.17)$$
But the right-hand side of the equality (3.17) is analytic in the right half -plane ${\rm Re}\ s > 0$.  This fact follows form the absolute convergence of the integral for ${\rm Re}\ s > 0$ and the uniform convergence with respect to $s,\  {\rm Re}\ s \ge x_0 > 0.$  Indeed,  since $\left(\mathcal{F} g\right) (t) \in L_q(\mathbb{R}),\ q= p/(p-1)$ we just apply the H$\ddot{o}$lder inequality to achieve the goal.  Hence,   the Euler integral for the gamma function and the Fubini theorem drive us at the equalities
$$\frac{4\  \Gamma(s) G_\alpha^*(s)}{ \Gamma\left((s + \alpha)/2\right)  \Gamma\left((s- \alpha)/2\right)}= 
  \int_{-\infty}^\infty  \left(\mathcal{F} g\right) (t) \int_0^\infty e^{- y\cosh t}y^{s-1} dy dt$$
$$=  \int_0^\infty  y^{s-1}  \int_{-\infty}^\infty  \left(\mathcal{F} g\right) (t) e^{- y\cosh t} dt dy.\eqno(3.18)$$
In the meantime, the Stirling asymptotic formula for the gamma function yields
$$\frac{ \Gamma(s) }{ \Gamma\left((s + \alpha)/2\right)  \Gamma\left((s- \alpha)/2\right)} =O\left( |s|^{1/2}\right),\  
|s| \to \infty.$$
Therefore under the condition $|s|^{1/2} G^*_\alpha(s) \in L_1(\mu-i\infty,\ \mu +i\infty),\ \mu >0$ the left-hand side of the first equality in (3.18) is integrable.  Taking the inverse Mellin transform (2.20), we find 
$$ {2\over \pi i} \int_{\mu-i\infty}^{\mu+i\infty} \frac{ \Gamma(s) G_\alpha^*(s)}{ \Gamma\left((s + \alpha)/2\right)  \Gamma\left((s- \alpha)/2\right)} y^{-s} ds =  \int_{-\infty}^\infty  \left(\mathcal{F} g\right) (t) e^{- y\cosh t} dt .\eqno(3.19)$$
Moreover, the right-hand side of (3.19) can be rewritten owing to the Parseval equality  for the Fourier transform \cite{tit} (see (2.12)) and integral representation (1.7) of the Macdonald function. Hence
$$ {1\over \pi i} \int_{\mu-i\infty}^{\mu+i\infty} \frac{ \Gamma(s) G_\alpha^*(s)}{ \Gamma\left((s + \alpha)/2\right)  \Gamma\left((s- \alpha)/2\right)} y^{-s} ds =  \int_{-\infty}^\infty  K_{i\tau}(y)  g(\tau) d\tau .\eqno(3.20)$$
The next step is to multiply both sides of (3.20) by
$$\frac{\Gamma(2\varepsilon) 2^{1-\varepsilon} }{\pi \Gamma(\varepsilon) \left|\Gamma(\varepsilon+ ix)\right|^2}
y^{\varepsilon-1} K_{ix}(y),\quad \varepsilon > 0, \  y >0, \   x \in \mathbb{R}$$
and integrate with respect to $y$ over $\mathbb{R}_+$.   Then changing the order of integration, recalling the Fubini theorem,  and calculating the inner integrals, appealing to relations (2.16.2.2) and (2.16.33.2) in \cite{prud}, Vol. II, we obtain
 $$ \frac{\Gamma(2\varepsilon) }{\pi^2 i  \Gamma(\varepsilon) \left|\Gamma(\varepsilon+ ix)\right|^2} \int_{\mu-i\infty}^{\mu+i\infty} \frac{ \Gamma(s) G_\alpha^*(s)}{ \Gamma\left((s + \alpha)/2\right)  \Gamma\left((s- \alpha)/2\right)} \Gamma\left({\varepsilon- s + ix\over 2} \right) \Gamma\left({\varepsilon- s -  ix\over 2} \right)    2^{-s-1} ds $$$$=  \frac{1 }{4 \pi B(\varepsilon, \varepsilon) }\int_{-\infty}^\infty   \left| B\left({\varepsilon + i(x- \tau)\over 2},\  {\varepsilon + i(x+ \tau)\over 2}\right) \right|^2  g(\tau) d\tau,\  0< \mu < \varepsilon .\eqno(3.21)$$
The left-hand side of  (3.21) can be expressed with the use of the Parseval equality (2.18).   In fact, we derive
$$ \frac{\Gamma(2\varepsilon) }{\pi^2 i  \Gamma(\varepsilon) \left|\Gamma(\varepsilon+ ix)\right|^2} \int_{\mu-i\infty}^{\mu+i\infty} \frac{ \Gamma(s) G_\alpha^*(s)}{ \Gamma\left((s + \alpha)/2\right)  \Gamma\left((s- \alpha)/2\right)} \Gamma\left({\varepsilon- s + ix\over 2} \right) \Gamma\left({\varepsilon- s -  ix\over 2} \right)    2^{-s-1} ds $$
$$=   \frac{\Gamma(2\varepsilon) }{\pi \Gamma(\varepsilon) \left|\Gamma(\varepsilon+ ix)\right|^2}
\int_0^\infty \hat{S}_{\alpha,\varepsilon} \left({1\over t},\ x\right) G_\alpha(t) {dt\over t},$$ 
where analogously to the calculation of the kernel $S_\alpha(x,\tau)$ (3.14)
$$\hat{S}_{\alpha,\varepsilon}\left(u,\ x\right) = \frac{1 }{2 \pi \sqrt \pi\  i  } \int_{\mu/2-i\infty}^{\mu/2 +i\infty} \frac{ \Gamma(s)\Gamma(1/2+s) }{ \Gamma\left(s + \alpha/2\right)  \Gamma\left(s- \alpha/2\right)} \Gamma\left({\varepsilon+ ix\over 2}- s  \right) \Gamma\left({\varepsilon -  ix\over 2} - s\right)    u^{-2s} ds$$
$$=  \frac{2 (2u)^{-\varepsilon-ix} \Gamma(\varepsilon+ix)\Gamma(-ix)}{  \Gamma\left((\varepsilon+\alpha+  ix)/2 \right) \Gamma\left((\varepsilon +  ix-\alpha)/2\right)}  {}_2F_3\left( {\varepsilon + ix\over 2},\    {\varepsilon+  ix + 1\over 2};  \ 1+ ix,    \   {\varepsilon+ \alpha+ ix\over  2},\    { \varepsilon+ ix - \alpha\over  2} ;  \   {1\over u^2} \right)$$
 $$ + 
 \frac{2 (2u)^{-\varepsilon+ix} \Gamma(\varepsilon- ix)\Gamma(ix)}{  
 \Gamma\left((\varepsilon+\alpha-  ix)/2 \right) \Gamma\left((\varepsilon -  ix-\alpha)/2\right)} {}_2F_3\left( {\varepsilon - ix\over 2},\    {\varepsilon-  ix + 1\over 2};  \ 1- ix,    \   {\varepsilon+ \alpha- ix\over  2},\    { \varepsilon- ix - \alpha\over  2} ;  \   {1\over u^2} \right).\eqno(3.22)$$

Further, returning to (3.21), we consider its right-hand side
$$I(\varepsilon, x) = \frac{1 }{4 \pi B(\varepsilon, \varepsilon) }\int_{-\infty}^\infty   \left| B\left({\varepsilon + i(x- \tau)\over 2},\  {\varepsilon + i(x+ \tau)\over 2}\right) \right|^2  g(\tau) d\tau.\eqno(3.23)$$
Since the square of the modulus of the beta function is even with respect to $\tau$ and $x$,  one can assume without loss of generality that $g(\tau)$ is an even function and $x >  0$.   Our goal is to prove the limit equality for some positive  $x$, where $g(x)$ satisfies the Dini condition,  namely, 
$$\lim_{\varepsilon \to 0+} I(\varepsilon, x) = g(x),\ x >0.\eqno(3.24)$$
In fact, appealing to (3.15),  we write for some small positive $\delta$
$$I(\varepsilon, x)-  g(x)=  \frac{1 }{4 \pi B(\varepsilon, \varepsilon) }\int_{-\infty}^\infty   \left| B\left({\varepsilon + i(x- \tau)\over 2},\  {\varepsilon + i(x+ \tau)\over 2}\right) \right|^2 \left[  g(\tau) - g(x)\right] d\tau$$
$$= \frac{4( (1+\varepsilon)^2+ x^2)( \varepsilon^2+ x^2) }{ \pi B(\varepsilon, \varepsilon) }\int_{-\infty}^\infty   \left| B\left(1+ {\varepsilon + i(x- \tau)\over 2},\  1+ {\varepsilon + i(x+ \tau)\over 2}\right) \right|^2 \frac{\left[  g(\tau) - g(x)\right]}
{(\varepsilon^2+ (x+\tau)^2) (\varepsilon^2+ (x-\tau)^2)} d\tau$$
$$= \frac{4( (1+\varepsilon)^2+ x^2)( \varepsilon^2+ x^2) }{ \pi B(\varepsilon, \varepsilon) }
\left(\int_{-\infty}^{-x-\delta}+ \int_{-x-\delta}^{-x+\delta} +  \int_{-x+ \delta}^{x- \delta} + \int_{x-\delta}^{x+\delta} + \int_{x+\delta}^ \infty \right)  \left| B\left(1+ {\varepsilon + i(x- \tau)\over 2},\  1+ {\varepsilon + i(x+ \tau)\over 2}\right) \right|^2$$
$$\times  \frac{\left[  g(\tau) - g(x)\right]} {(\varepsilon^2+ (x+\tau)^2) (\varepsilon^2+ (x-\tau)^2)} d\tau
= I_1(\varepsilon, x)+  I_2(\varepsilon, x)+  I_3(\varepsilon, x)+   I_4(\varepsilon, x)+   I_5(\varepsilon, x). $$

Starting from integral $I_5(\varepsilon, x)$,  we find
$$|I_5(\varepsilon, x)|  \le   \frac{4B^2 (1+ \varepsilon/2,\  1+\varepsilon/2) ( (1+\varepsilon)^2+ x^2)( \varepsilon^2+ x^2) }{ \pi B(\varepsilon, \varepsilon) } \int_{x+\delta}^ \infty   \frac{\left| g(\tau) - g(x)\right|} {(\varepsilon^2+ (x+\tau)^2) (\varepsilon^2+ (x-\tau)^2)} d\tau $$
$$=  \frac{2^{2\varepsilon+1} B^2 (1+ \varepsilon/2,\  1+\varepsilon/2) ( (1+\varepsilon)^2+ x^2)( \varepsilon^2+ x^2)\Gamma(\varepsilon + 1/2) }{ \pi \sqrt \pi \Gamma(\varepsilon) } \int_{\delta}^ \infty   \frac{\left| g(x+t) - g(x)\right|} {(\varepsilon^2+ (2x+ t)^2) (\varepsilon^2+ t^2)} dt $$
$$\le C_1(x)  \varepsilon   \int_{\delta}^ \infty   \frac{| g(x+t)| +|g(x)|} {t^4} dt \le  C_2(x)  {\varepsilon \over \delta^3} + 
C_1(x)\varepsilon \left( \int_{ \delta}^ \infty  t^{-4q} dt \right)^{1/q}  ||g||_{L_p(\mathbb{R})}$$
$$=  {\varepsilon\over \delta^3}  \  \left( C_2(x) +  {C_3(x) \over \delta^{1- 1/q } }\right) ,$$
where $C_i(x),\ i=1,2,3$ are constants.  Clearly, one can make the latter expression arbitrary small, choosing first some $\delta$ and then $\varepsilon$.  Thus 
$$\lim_{\varepsilon \to 0+} I_5(\varepsilon, x) = 0.$$
In the same manner we establish the equality 
$$\lim_{\varepsilon \to 0+} I_1(\varepsilon, x) = 0.$$
Concerning integral $I_3(\varepsilon, x)$,  we presume that $\delta $ is small and does not exceed $x$.    Hence, analogously, 
$$|I_3(\varepsilon, x)|   \le \varepsilon \  C_4(x)   \int_{\delta}^{2x- \delta}  \frac{| g(x-t)| +|g(x)|} {t^2 (2x-t)^2} dt 
\le {\varepsilon\over \delta^3}  \  \left[ C_5(x)+ {C_6(x)\over \delta^{1-1/q} }\right]    $$
and again 
$$\lim_{\varepsilon \to 0+} I_3(\varepsilon, x) = 0.$$
Finally, we estimate integrals $I_2(\varepsilon, x),\  I_4(\varepsilon, x)$.  We have 
$$|I_4(\varepsilon, x)|\le \varepsilon \  C_7(x) \int_{x-\delta}^{x+\delta}   \frac{\left| g(\tau) - g(x)\right|} {(\varepsilon^2+ (x+\tau)^2) (\varepsilon^2+ (x-\tau)^2)} d\tau =  \varepsilon \  C_7(x) \int_{-\delta}^{\delta}   \frac{\left| g(x+t) - g(x)\right|} {(\varepsilon^2+ (2x+t)^2) (\varepsilon^2+ t^2)} dt$$
$$\le {\varepsilon\over x^2}  \  C_7(x) \int_{-\delta}^{\delta}   \frac{\left| g(x+t) - g(x)\right|} {\varepsilon^2+ t^2} dt
\le {1\over 2 x^2}  \  C_7(x) \int_{-\delta}^{\delta}   \frac{\left| g(x+t) - g(x)\right|} {| t|} dt,$$
where the latter integral becomes small when $\delta$  goes to zero via the absolute continuity of the Lebesgue integral. 
Therefore $I_4(\varepsilon, x) \to 0, \  \varepsilon \to 0+$ and in the same manner we get 
$$\lim_{\varepsilon \to 0+} I_2(\varepsilon, x) = 0.$$
Thus (3.24) holds true,  and the inversion formula (3.16) follows immediately after the passage to the limit in (3.21) by $\varepsilon \to 0+$.    

\end{proof} 

\section{Particular cases}

\subsection{The case $\alpha =0$.}     Letting $\alpha =0$, we come up with the Lebedev type index transform, involving  a square of the Macdonald function
$$F_0(\tau) = \int_0^\infty  K^2_{i\tau/2}(x)  f(x)dx,   \quad    \tau \in \mathbb{R}. \eqno(4.1)$$
In this form the transform is mentioned by formula (8.59) in \cite{mar} as a particular case of the general Wimp-Yakubovich transform with respect an  index of the Meijer $G$-function (cf.  \cite{yal},  Chapter 7).    Despite Theorem 6 is proved for non-zero $\alpha$,  one can adjust its proof for the zero- case too.  In fact, the hypergeometric function in the inversion formula (3.1) is reduced to the ${}_1F_2$- function  and can be expressed in terms of the modified Bessel functions.  Precisely, appealing to relation (7.14.1.4) in \cite{prud}, Vol. III, we find 
$${}_1F_2\left(  {1+ i\tau\over 2};  \ 1+ i\tau,    \   { i\tau\over  2} ;  \   x^2 \right) = 2 \ I_{i\tau/2} (x) \Gamma\left(1+ {i\tau\over 2}\right) \left({x\over 2}\right)^{-i\tau}  \  \left[ {x\over 2} \  \Gamma\left({i\tau\over 2}\right)  I_{i\tau/2  -1}(x)\right.$$
$$\left.  -   \Gamma\left(1+ {i\tau\over 2}\right) I_{i\tau/2} (x) \right].$$
Hence
$$\frac{   (x/2)^{i\tau-1}}{\Gamma^2\left(i\tau/ 2\right) }  \  {}_1F_2\left( {1+ i\tau\over 2};  \ 1+ i\tau,  \   { i\tau\over  2} ;  \   x^2 \right) = i\tau \   I_{i\tau/2} (x) \left[  I_{i\tau/2  -1}(x) - {i\tau \over 2x}  \    I_{i\tau/2} (x) \right].\eqno(4.2)$$ 
Meanwhile  calling properties of the modified Bessel functions (see in \cite{erd},  Vol. II), we obtain
$$ {i\tau \over x}  \    I_{i\tau/2} (x) =  I_{i\tau/2  -1}(x) -   I_{i\tau/2  +1}(x),$$
$$ I_{i\tau/2  -1}(x) +  I_{i\tau/2  +1}(x) = 2 {d\over dx} I_{i\tau/2}(x).$$ 
Substituting these expressions in (4.2),  after straightforward  simplifications we derive
$$\frac{   (x/2)^{i\tau-1}}{\Gamma^2\left(i\tau/ 2\right) }  \  {}_1F_2\left( {1+ i\tau\over 2};  \ 1+ i\tau,  \   { i\tau\over  2} ;  \   x^2 \right) =   i\tau \   I_{i\tau/2} (x)  {d\over dx} I_{i\tau/2}(x) = {i\tau\over 2}   {d\over dx} I^2_{i\tau/2}(x),$$ 
and the inversion formula (3.1) for the index transform (4.1) takes the form
$$f(x) =  { i \over  2\pi  }   \int_{-\infty} ^\infty  \   {d\over dx} I^2_{i\tau/2}(x)\      F_0 (\tau)\  \tau \ d\tau. \eqno(4.3)$$
However this is exactly the inversion formula (8.60) in \cite{mar} subject to elementary changes  of variables  and functions.   An analog of Theorem 6 is 

{\bf Theorem 8}.  {\it  Let $   \nu < 1,\  1< p \le 2,\  q= p/(p-1). $    Let $s f^*(s) \in L_p(\nu-i\infty, \nu+i\infty)$, where $f^*(s)$ is the Mellin transform $(2.19)$ of $f \in L_1((1,\infty); \  dt)$. If, besides above assumptions,  the index  transform $(4.1)$  satisfies the condition $\tau e^{\pi |\tau|}  F_0(\tau) \in L_1(\mathbb{R})$ and the  Mellin transform $f^*(s)$ vanishes at the point $s=1$, then for all $x> 0$  inversion formula $(4.3)$ holds,  where the corresponding integral is absolutely convergent.}

\begin{proof}  The scheme of the proof is the same  as in Theorem 6.  Nevertheless, to adjust the value $\alpha=0$  and to have a strip,  where one can choose a vertical line for the integration in the right-hand side of (3.11),  it   can be rewritten as 
$${1\over \pi^3 i}  \int_{\mu- i\infty}^{\mu+i\infty} 2^{-s-2} \Gamma(s)   \int_{-\infty} ^\infty \tau \sinh(\pi\tau) \   F_0(\tau) \Gamma\left({1-s+ i\tau\over 2}\right)  \Gamma\left({1-s- i\tau\over 2}\right)  d\tau  \    \lambda^{-s} ds $$$$=   {1\over 2\pi i}  \int_{\mu- i\infty}^{\mu +i\infty}  \Gamma\left({s+1\over 2}\right)  \Gamma\left({s-1\over 2}\right)   f^*(2-s) \ \lambda^{-s} ds$$
with  $0 < \mu < 1$ owing to the analyticity, integrability conditions and the value $f^*(1)=0$.  Hence one can proceed all further steps of the proof of  Theorem 6.
\end{proof} 

The corresponding adjoint operator to (4.1) has the form (see (1.2))
$$G_0(x) = \int_0^\infty  K^2_{i\tau/2}(x)  g(\tau)d\tau. \eqno(4.4)$$

Hence Theorem 7 becomes

{\bf Theorem 9.} {\it Let $  g(\tau) \in L_p(\mathbb{R}),\ 1 < p \le 2,$  satisfying the Dini condition at some point $x \in \mathbb{R}\backslash \{0\} $. Let  also the index transform $(4.4)$ be such that its Mellin transform $ G^*_0(s) \in L_1((\mu-i\infty,\  \mu +i\infty);  |s|^{1/2} ds ),\   \mu >0$.  Then 
$$ g(x) = \lim_{\varepsilon \to 0+} {1\over \pi} \  \frac{ \Gamma(2\varepsilon) }{ \Gamma(\varepsilon)}
\int_0^\infty \left[\frac{  \Gamma(-ix)\left(t/2\right) ^{\varepsilon +ix-1} }{ \Gamma(\varepsilon- ix) \Gamma^2\left((\varepsilon+  ix)/2 \right)} \   {}_1F_2\left(  {\varepsilon+  ix + 1\over 2};  \ 1+ ix,    \   {\varepsilon+  ix\over  2} ;  \   t^2 \right)\right.$$
$$\left.    + \frac{  \Gamma(ix)\left(t/ 2\right)^{\varepsilon - ix-1} }{ \Gamma(\varepsilon+ ix) \Gamma^2\left((\varepsilon-  ix)/2 \right) }  \   {}_1F_2\left(  {\varepsilon-  ix + 1\over 2};  \ 1- ix,    \   {\varepsilon- ix\over  2} ;  \   t^2 \right) \right]G_0(t) dt, \eqno(4.5)$$
where the limit is pointwise.} 

{\bf Remark 1}.  When it is possible to pass to the limit under the integral sign in (4.5), this formula takes the form
$$g(x)= {ix\over 2\pi} \int_0^\infty  {d\over dt }\left[ I^2_{ix/2}(t)-  I^2_{- ix/2}(t)\right] G_0(t) dt.$$
Moreover,  with elementary changes of the variables,  functions as in the previous example and the integration by parts we arrive at the pair of Lebedev transforms  \cite{square}  (see formulas (8.57), (8.58) in \cite{mar}). 

\subsection{The case $\alpha=1$.}   In this case we will consider the following index transforms 

$$F_1(\tau) = \int_0^\infty \left| K_{(i\tau+ 1)/2}(x)\right|^2  f(x)dx,   \quad    \tau \in \mathbb{R}, \eqno(4.6)$$
$$G_1(x) = \int_{-\infty}^\infty \left| K_{(i\tau+ 1)/2}(x)\right|^2  g(\tau)d\tau,   \quad   x  \in \mathbb{R}_+. \eqno(4.7)$$
Employing relation (7.14.1.5) in \cite{prud}, Vol. III, the corresponding kernel in the inversion formula (3.1) can be reduced to the product of the modified Bessel functions. Precisely, we find
$$\frac{   (x/2)^{i\tau-1}}{\Gamma\left((i\tau- 1)/ 2\right)  \Gamma\left((1+ i\tau)/ 2\right)}  \  {}_1F_2\left( {i\tau\over 2};  \ 1+ i\tau,   {i\tau-1 \over 2} ;  \   x^2 \right) =   I_{(i\tau-1)/2} (x)$$
$$\times \left[  {i\tau-1\over 2}   I_{(i\tau-1)/2} (x) + {x\over 2} I_{(i\tau+1)/2}(x)\right] = {x\over 2} \   I_{(i\tau-1)/2} (x) I_{(i\tau-3)/2} (x).$$

Thus the inversion theorems for operators (4.6), (4.7) are given accordingly,

{\bf Theorem 10}.  {\it  Let $ \nu <  0,\  1< p \le 2,\  q= p/(p-1). $    Let $s f^*(s) \in L_p(\nu-i\infty, \nu+i\infty)$, where $f^*(s)$ is the Mellin transform $(2.19)$ of $f \in L_1((1,\infty); \  t dt)$.  If, besides above assumptions,  the index transform  $(4.6)$ of $f$ satisfies the condition $\tau e^{\pi |\tau|}  F_1(\tau) \in L_1(\mathbb{R})$ and  $f^*(0)=0$, then for all $x> 0$ the following inversion formula holds
$$f(x) =  { 1\over  \pi  }   \int_{-\infty} ^\infty \left[  {x\over 2} \   I_{(i\tau-1)/2} (x) I_{(i\tau-3)/2} (x)   +  \   \frac{\cosh(\pi\tau/2) }{\pi   x }\right] \      F_1(\tau)\  d\tau,$$
and the integral converges absolutely.}

{\bf Theorem 11.} {\it Let $ g(\tau) \in L_p(\mathbb{R}),\ 1 < p \le 2,$  satisfying the Dini condition at some point $x \in \mathbb{R}\backslash \{0\} $. Let  also the index transform $G_1$ be such that its Mellin transform $ G^*_1(s) \in L_1((\mu-i\infty,\  \mu +i\infty);  |s|^{1/2} ds ),\   \mu >0$.  Then 
$$ g(x) = \lim_{\varepsilon \to 0+} {1\over \pi} \  \frac{ \Gamma(2\varepsilon) }{ \Gamma(\varepsilon)}
\int_0^\infty \left[\frac{  \Gamma(-ix)\left(t/2\right) ^{\varepsilon +ix-1} }{ \Gamma(\varepsilon- ix) \Gamma\left((\varepsilon+1+  ix)/2 \right) \Gamma\left((\varepsilon +  ix- 1)/2\right)} \right.$$$$\left. \times  {}_1F_2\left( {\varepsilon + ix\over 2}    ;  \ 1+ ix,    \     { \varepsilon+ ix - 1\over  2} ;  \   t^2 \right)   + \frac{  \Gamma(ix)\left(t/ 2\right)^{\varepsilon - ix-1} }{ \Gamma(\varepsilon+ ix) \Gamma\left((\varepsilon+ 1-  ix)/2 \right) \Gamma\left((\varepsilon -  ix- 1)/2\right)} \right.$$$$\left. \times  {}_1F_2\left( {\varepsilon - ix\over 2};  \ 1- ix,    \   { \varepsilon- ix - 1\over  2} ;  \   t^2 \right) \right]G_1(t) dt,\eqno(4.8) $$
where the limit is pointwise.  When the passage to the limit under the integral sign is possible, the inversion formula $(4.8)$ takes the form} 
$$ g(x) =  {1\over 4\pi} \int_0^\infty \left[  I_{(ix -1)/2} (t) I_{(ix -3)/2} (t)+   I_{- (ix+ 1)/2} (t) I_{-(i\tau+3)/2} (t)\right]G_1(t) \ t dt. $$

\section{Initial   value problem}

The index transform (1.2) can be successfully applied to solve a  boundary problem  for the following partial differential difference equation, involving the Laplacian
$$\Delta u_n= u_{n+2}+ 2 u_n + u_{n-2}, \eqno(5.1)$$ 
where $\Delta = {\partial^2 \over \partial x^2} +  {\partial^2 \over \partial y ^2}$ is the Laplacian in $\mathbb{R}^2,\ n \in \mathbb{Z}$  and $u_n= u_{-n}$.   Concerning properties and solutions of the ordinary differential difference equations see the survey \cite{bat}. Writing equation (5.1) in polar coordinates $(r,\theta)$, precisely
$$ {\partial^2 u_n \over \partial r^2}+  {1\over r}   {\partial u_n \over \partial r}+  {1\over r^2}   {\partial^2 u_n \over \partial \theta^2} =  u_{n+2}+ 2 u_n + u_{n-2}, \eqno(5.2)$$
we establish the following 

{\bf Lemma 4.} {\it  Let $g \in L_1\left(\mathbb{R}; e^{ (2\pi- \beta)|\tau|} d\tau\right),\  \beta \in [0, \pi/2[$. Then  functions
 $$u_n(r,\theta) =  \int_{-\infty}^\infty e^{\theta\tau} \left| K_{(i\tau+ n)/2}(r)\right|^2  g(\tau)d\tau,\quad  n \in \mathbb{Z},\eqno(5.3)$$
where $r >0, \ 0\le \theta \le 2\pi$ satisfy  the partial  differential difference equation $(5.2)$, vanishing at infinity.}

\begin{proof} In fact, this follows from the direct substitution (5.3) into (5.2) and the use of Lemma 2.  The necessary  differentiation  with respect to $r$ under integral sign is allowed via the absolute and uniform convergence, which can be justified  employing the inequality \cite{yak}  $(z=\mu+i\tau)$
$$\left|K_z(x)\right| \le e^{-\beta |\tau|} K_\mu(x\cos\beta),\   x >0, \   \beta \in [0, \pi/2[$$
and the asymptotic behavior  (1.4) at infinity  of the Macdonald function.  
\end{proof}

Finally, as a direct consequence of Theorem 7, we will formulate the initial value problem for equation (5.2) and give its solutions.

{\bf Theorem 12.} {\it Let $n \in \mathbb{Z},\  G_n(r)$ and 
$$g(x) =  \lim_{\varepsilon \to 0+} {1\over \pi} \  \frac{ \Gamma(2\varepsilon) }{ \Gamma(\varepsilon)}
\int_0^\infty \left[\frac{  \Gamma(-ix)\left(t/2\right) ^{\varepsilon +ix-1} }{ \Gamma(\varepsilon- ix) \Gamma\left((\varepsilon+ n+  ix)/2 \right) \Gamma\left((\varepsilon +  ix- n)/2\right)} \right.$$$$\left. \times  {}_2F_3\left( {\varepsilon + ix\over 2},\    {\varepsilon+  ix + 1\over 2};  \ 1+ ix,    \   {\varepsilon+ n + ix\over  2},\    { \varepsilon+ ix - n \over  2} ;  \   t^2 \right)   + \frac{  \Gamma(ix)\left(t/ 2\right)^{\varepsilon - ix-1} }{ \Gamma(\varepsilon+ ix) \Gamma\left((\varepsilon+ n-  ix)/2 \right) \Gamma\left((\varepsilon -  ix- n)/2\right)} \right.$$$$\left. \times  {}_2F_3\left( {\varepsilon - ix\over 2},\    {\varepsilon-  ix + 1\over 2};  \ 1- ix,    \   {\varepsilon+ n- ix\over  2},\    { \varepsilon- ix - n\over  2} ;  \   t^2 \right) \right]G_n(t) dt$$
satisfy conditions of Theorem $7$. Then  functions $u_n(r,\theta),\   r >0,  \  0\le \theta < \pi/2$ by formula $(5.3)$  will be solutions  of the initial value problem for the partial differential difference equation $(5.2)$ subject to the initial condition}
$$u_n(r,0) = G_n (r).$$

\bigskip
\centerline{{\bf Acknowledgments}}
\bigskip
The present investigation was supported, in part,  by the "Centro
de Matem{\'a}tica" of the University of Porto.

\bibliographystyle{amsplain}

\begin{thebibliography}{10}

\bibitem{yak}  S. Yakubovich, {\it Index Transforms},  World Scientific Publishing Company, Singapore, New Jersey, London and Hong Kong (1996).

\bibitem{yal}  S. Yakubovich and Yu.  Luchko, {\it The Hypergeometric Approach to Integral Transforms and Convolutions}, (Kluwers Ser.  Math. and Appl.: Vol. 287), Dordrecht, Boston, London (1994).


\bibitem{erd}   A. Erd\'elyi, W. Magnus, F. Oberhettinger and F.G. Tricomi, {\it Higher Transcendental Functions}, Vols. I,  II, McGraw-Hill, New  York, London and Toronto (1953).

\bibitem{prud} A.P. Prudnikov, Yu.A. Brychkov and O.I. Marichev, \textit{Integrals and Series:  Vol. I: Elementary Functions}, Gordon and Breach, New York,  1986;  {\it Vol. II:  Special Functions}, Gordon and Breach, New York  (1986);  {\it Vol. III:  More Special Functions},   Gordon and Breach, New York  (1990).

\bibitem{square} N.N. Lebedev,   On an integral representation of an arbitrary function in terms of squares of Macdonald functions with imaginary index,  {\it Sibirsk. Mat. Zh.},   {\bf 3}  (1962),  213- 222 (in Russian).

\bibitem{leb} S. Yakubovich,  On Lebedev type integral transformations associated with modified Bessel functions, {\it Nieuw Arch. Wisk.},  (4) {\bf 17}  (1999), no. 2,  219-  227.

\bibitem{yaprod} S. Yakubovich,  On a new index transformation related to the product of Macdonald functions,  {\it Rad. Mat.},  {\bf 13}  (2004), no. 1, 63 - 85. 

\bibitem{yaprodrev} S. Yakubovich,   A double index transform with a product of Macdonald's functions revisited,  {\it Opuscula Math.}  {\bf 29}  (2009), no. 3,  31- 329. 


\bibitem{yakl} S. Yakubovich,  On the Kontorovich-Lebedev transformation, \textit{J. Integral Equations Appl.}  {\bf 15}  (2003), no. 1, 95 - 112.


\bibitem{gen} S. Yakubovich,  $L_p$ -boundedness of general index transforms,  {\it Liet. Mat. Rink.},  45 (2005), no. 1, 127- 147; translation in Lithuanian Math. J. {\bf 45}  (2005), no. 1, 102-  122.

\bibitem {tit}  E.C. Titchmarsh, {\it  An Introduction to the Theory of Fourier Integrals},    Chelsea, New York  ( 1986).

\bibitem {mar}  O.I. Marichev, {\it  Handbook of Integral Transforms of Higher Transcendental Functions. Theory and Algorithmic Tables},    Chichester: Ellis Horwood   ( 1983).

\bibitem{bat} H. Bateman,  Some simple differential difference equations and the related functions,  {\it Bull. Amer. Math. Soc.},  49 (1943), 494- 512.





\end{thebibliography}

\end{document}